\mag=\magstephalf
\input amstex
\documentstyle{amsppt}
\pagewidth{12.6cm}
\pageheight{19.8cm}
\NoRunningHeads
\hcorrection{1cm}
\def\t#1{\text{\rm #1}}
\def\thetag#1{(#1)}
\def\vep{\varepsilon}
\def\sig{\sigma}
\def\alp{\alpha}

\def\lam{\lambda}
\def\Lam{\Lambda}

\def\gam{\gamma}
\def\vth{\vartheta}

\def\lot{\overset . \to\otimes}
\def\hotimes{\overset . \to{\circ}}

\def\T{\Tilde}
\def\tr{\operatorname{tr}}

\def\Ch{\operatorname{Ch}}
\def\dim{\operatorname{dim}}
\def\vdim{ \operatorname {\bold{dim}}}
\def\diag{ \operatorname {diag}}
\def\Irr{ \operatorname {Irr}}
\def\id{\operatorname {id}}
\def\Id{\operatorname {Id}}
\def\E{ \operatorname {End}}
\def\H{ \operatorname {Hom}}

\def\lfl{\lfloor}
\def\rfl{\rfloor}

\def\fS{\frak S}
\def\fq{\frak q}
\def\fgl{\frak g \frak l}

\def\fq{{\frak q}}

\def\pr{{\frak p}}

\def\A{{\Cal A}}
\def\B{{\Cal B}}
\def\C{{\Cal C}}
\def\U{{\Cal U}}

\topmatter
\title 
A duality of a twisted group algebra of the hyperoctahedral group 
and the queer Lie superalgebra 
\endtitle
\author
Manabu Yamaguchi
\endauthor
\affil
Department of Mathematics, Aoyama Gakuin University\\
Chitosedai 6-16-1, Setagaya-ku, Tokyo, 157-8572 Japan 
\endaffil
\endtopmatter
\document						
\head
\S1. Introduction
\endhead
We establish a duality relation (Theorem 4.2) between 
one of the twisted group algebras, 
of the hyperoctahedral group $H_k$ (or the Weyl group of type $B_k$) and a Lie 
superalgebra $\fq(n_0)\oplus\fq(n_1)$ for any integers $k\geq4$ and 
$n_0$, $n_1\geq1$. 
Here $\fq(n_0)$ and $\fq(n_1)$ denote the ``queer'' Liesuperalgebras as called 
by some authors. The twisted group algebra $\B'_k$ in focus in this paper 
belongs to a different cocycle from the one $\B_k$ used by 
A\. N\. Sergeev in his work [8] on a duality with $\fq(n)$ 
and by the present author in a previous work 
[11]. This $\B'_k$ contains the twisted group algebra $\A_k$ 
of the symmetric group $\fS_k$ in a straightforward manner (\S1\. 1\. 1), and 
has a structure similar to the semidirect product of $\A_k$ and 
${\Bbb C}[(\Bbb Z/2{\Bbb Z})^k]$. ($\B'_k$ and $\B_k$ were denoted by 
${\Bbb C}^{[-1,+1,+1]}W_k$ and 
${\Bbb C}^{[+1,+1,-1]}W_k$ respectively by J\. R\. Stembridge in [10].) 

In \S2, we construct the ${\Bbb Z}_2$-graded simple $\B'_k$-modules (where 
${\Bbb Z}_2={\Bbb Z}/2{\Bbb Z}$) using 
an analogue of the little group method. These simple $\B'_k$ modules 
are slightly different from the non-graded simple $\B'_k$-modules 
constructed by Stembridge in [10] because of the difference between 
${\Bbb Z}_2$-graded and non-graded theories, but they can easily be translated 
into each other. We will use the algebra $\C_k\lot\B'_k$, where 
$\C_k$ is the $2^k$-dimensional Clifford algebra (cf\. \thetag{3.2}) 
and $\lot$ denotes the ${\Bbb Z}_2$-graded tensor product 
(cf\. [1], [2], [11, \S1]), as an intermediary for establishing our duality, 
as we explain below. The construction of the simple $\B'_k$-modules leads to 
a construction of the simple $\C_k\lot\B'_k$-modules in \S3. 

In \S4, we define a representation of $\C_k\lot\B'_k$ in 
the $k$-fold tensor product $W=V^{\otimes k}$ 
of $V={\Bbb C}^{n_0+n_1}\oplus{\Bbb C}^{n_0+n_1}$, 
the space of the natural representation 
of the Lie superalgebra $\fq(n_0+n_1)$. 
This representation of $\C_k\lot\B'_k$ depends on $n_0$ and $n_1$, not just 
$n_0+n_1$. 
\comment
the decomposition of $V$ into a direct sum of two graded spaces 
$X$ and $Y$, namely $V=X\oplus Y$, such that $\vdim X=(n_0,n_0)$, 
$\vdim Y=(n_1,n_1)$, and $n_0+n_1=n$. This representation is an extension 
of the representation of $\B_k$ in $W$ which was defined by A\. N\. Sergeev 
(cf\. Theorem 4.6). 
\endcomment
Note that $\B_k$ can be regarded as a subalgebra of $\C_k\lot\B'_k$, 
since $\B_k$ is isomorphic to $\C_k\lot \A_k$ by our previous result 
(cf\. \thetag{3.3} of [11]). Under this embedding, our representation 
of $\C_k\lot\B'_k$ restricts to the representation of $\B_k$ in $W$ 
defined by Sergeev (cf\. Theorem A). 
We show that the centralizer of $\C_k\lot\B'_k$ in $\E(W)$ 
is generated by the action of the Lie superalgebra 
$\fq(n_0)\oplus\fq(n_1)$ (Theorem 4.1). 
Moreover we show that 
$\B'_k$ and $\fq(n_0)\oplus\fq(n_1)$ 
act on a subspace $W'$ of $W$ ``as mutual centralizers of each other'' 
(Theorem 4.2). 
Note that $\A_k$ and $\fq(n)$ 
act on the same space $W'$ 
``as mutual centralizers of each other'' (cf\. Theorem B).

In Appendix, we include short explanations of some known reults, 
which we use in the previous sections. 

In this paper, all vector spaces, and associative algebras, 
and representations are assumed to 
be finite dimensional over ${\Bbb C}$ unless specified otherwise. 
The precise statements of the results skeched in the introduction 
use the formulation of ${\Bbb Z}_2$-graded representations of 
${\Bbb Z}_2$-graded algebras (superalgebras) 
(cf\. \S1.1.3) as was used in [1] and [2].

\subhead
1.1 Preliminaries
\endsubhead
\subsubhead
1.1.1. A twisted group algebra $\B'_k$
\endsubsubhead
For any $k\geq1$, let $\B'_k$ denote the associative algebra generated by 
$\tau'$ and the $\gamma_i$, $1\leq i\leq k-1$, with relations 
$$
\align
&{\tau'}^2=\gamma_i^2=1\quad(1\leq i\leq k-1),
\quad(\gamma_i\gamma_{i+1})^3=1\quad(1\leq i\leq k-2),\tag1.1\\
&(\gamma_i\gamma_j)^2=-1\quad(|i-j|\geq2),\quad
(\tau'\gamma_i)^2=1\quad(2\leq i\leq k-1),\\
&(\tau'\gamma_1)^4=1.\\
\endalign
$$ 
If $k\geq4$, then $\B'_k$ is isomorphic to 
a twisted group algebra of the hyperoctahedral group $H_k$ with a 
non-trivial $2$-cocycle (cf\. [10, Prop\. 1.1]). 
We regard $\B'_k$ as a ${\Bbb Z}_2$-graded algebra by giving 
the generator $\tau'$ (resp\. the generator $\gamma_i$, $1\leq i\leq k-1$) 
degree $0$ (resp\. degree $1$). 
Note that this grading of $\B'_k$ is different 
from that of $\B_k$ in \thetag{3.1} or in [11]. 

Let $\A_k$ denote the ${\Bbb Z}_2$-graded subalgebra of $\B'_k$ 
generated by $\gamma_i$, $1\leq i\leq k-1$. If $k\geq 4$, then $\A_k$ is 
isomorphic to a twisted group algebra of the symmetric group 
$\fS_k$ with a non-trivial $2$-cocycle, with the ${\Bbb Z}_2$-grading as in [2] 
and [11]. 
\subsubhead
1.1.2. Partitions and symmetric functions
\endsubsubhead
Let $P_k$ denote the set of all partitions of $k$, and 
put $P=\coprod_{k\geq0}P_k$. 
For $\lam\in P$, we write $l(\lam)$ for the length of $\lam$, namely 
the number of non-zero parts of $\lam$. Also we write $|\lam|=k$ 
if $\lam\in P_k$.
Let $DP_k$ and $OP_k$ denote the distinct partitions 
(or strict partitions, namely partitions whose parts are distinct) 
and the odd partitions (namely partitions whose parts are all odd) 
of $k$ respectively. 
Let $DP^{+}_k$ and $DP^{-}_k$ be the sets of all 
$\lam\in DP_k$ 
such that $(-1)^{k-l(\lam)}=+1$ and $-1$ respectively. Note that 
$(-1)^{k-l(\lam)}$ equals the signature of permutations with 
cycle type $\lam$. 
We also put $DP=\coprod_{k\geq0}DP_k$ and $OP=\coprod_{k\geq0}OP_k$. 
Let $(DP^2)_k$ (resp\. ${(OP^2)}_k$) denote the set 
of all $(\lam,\mu)\in DP^2$ (resp\. $OP^2$) 
such that $|\lam|+|\mu|=k$. 
Let ${(DP^2)}^{+}_k$ and ${(DP^2)}^{-}_k$ be the sets of all 
$(\lam,\mu)\in (DP^2)_k$ such that 
$(-1)^{k-l(\lam)-l(\mu)}=+1$ and $-1$ respectively. 

Let $\Lam_x$ denote the ring of the symmetric functions 
in infinitely many variables 
$x=\{x_1,x_2,\dots\}$ with coefficients in ${\Bbb C}$; 
namely our $\Lam_x$ is the scalar extension of the $\Lam_x$ in [6], which is 
${\Bbb Z}$-algebra, to ${\Bbb C}$. 

Let $\Omega_x$ denote the subring of 
$\Lambda_x$ generated by the power sums of odd degrees, namely the $p_r(x)$, 
$r=1,3,5$, $\dots$ . 
Then $\{p_{\lam}(x)\,|\,\lam\in OP\}$ is a basis of $\Omega_x$, where 
$p_{\lam}=\prod_{i\geq1}p_{\lam_i}$. 
For $\lam\in DP$, let $Q_{\lam}(x)\in \Lam_x$ 
denote Schur's $Q$-function indexed by $\lam$ (cf\. [7], [9, \S6]). 
Then $\{Q_{\lam}(x)\;|\;\lam\in DP\}$ 
is also a basis of $\Omega_x$. 
\subsubhead
1.1.3. Semisimple superalgebras
\endsubsubhead
This theory was developed by T\. J\'ozefiak in [1], which we mostly follow. 
A ${\Bbb Z}_2$-graded algebra $A$, which is called 
a {\bf superalgebra} in this paper, 
is called {\bf simple} if it does not have non-trivial ${\Bbb Z}_2$-graded 
two-sided ideals. 
If $A$ is a simple superalgebra, then it is either isomorphic to $M(m,n)$ 
(denoted by $M(m|n)$ in [2]) for some $m$ and $n$, 
or isomorphic to $Q(n)$ for some $n$ (see [2], [13, \S1] for 
the definitions of simple superalgebras $M(m,n)$, $Q(n)$). 

Let $V$ be an $A$-{\bf module}, namely a ${\Bbb Z}_2$-graded 
vector space $V=V_0\oplus V_1$ together with a representation 
$\rho\:A\to\E(V)$ satisfying $\rho(A_{\alp})V_{\beta}\subset V_{\alp+\beta}$ 
($\alp$, $\beta\in{\Bbb Z}_2$). 
We simply write $\rho(a)v=av$ for $a\in A$ and $v\in V$. 
By an $A$-submodule of $V$ we mean a ${\Bbb Z}_2$-graded $\rho(A)$-stable 
subspace of $V$. We say that $V$ is {\bf simple} if 
it does not have non-trivial $A$-submodules. 

Let $V$ and $W$ be two $A$-modules. Let 
$\H^{\alp}_A(V,W)$ $(\alp\in{\Bbb Z}_2)$ denote the subspace of 
$\H^{\alp}(V,W)=\{f\in\H(V,W)\,;\,f(V_{\beta})\subset W_{\alp+\beta}\}$ 
consisting of all elements $f\in \H^{\alp}(V,W)$ such that 
$f(av)=(-1)^{\alp\cdot\beta}af(v)$ 
for $a\in A_{\beta}$ $(\beta\in{\Bbb Z}_2)$, $v\in V$. 
Put $\H^{\cdot}_A(V,W)=\H^{0}_A(V,W)\oplus\H^{1}_A(V,W)$ and put 
$\E^{\cdot}_A(V)=\H^{\cdot}_A(V,V)$. 
We call $\E^{\cdot}_A(V)$ the {\bf supercentralizer } of 
$A$ in $\E(V)$. 
Two $A$-modules $V$ and $W$ are called {\bf isomorphic} 
if there exists an invertible linear map $f\in\H^{\cdot}_A(V,W)$. 
If this is the case, we write $V\cong_AW$ (or simply write $V\cong W$). 
If $V$ and $W$ are simple $A$-modules, then 
$V\cong W$ if and only if there exists 
an invertible element in $\H^{0}_A(V,W)$ or $\H^{1}_A(V,W)$. 
Note that, in [11] we distinguished between $V$ and 
the shift of $V$ which is defined 
to be the same vector space as $V$ with the switched grading. 
In this paper, however, we identify $V$ and the shift of $V$.

If $V$ is a simple $A$-module, 
then $\E^{\cdot}_A(V)$ is isomorphic to either 
$M(1,0)\cong{\Bbb C}$ or $Q(1)\cong\C_1$ (cf\. [1, Prop\. 2.17], 
[2, Prop\.  2.5, Cor\.  2.6]). 
In the former (resp\. latter) case, we say 
that $V$ is of {\bf type} $M$ (resp\. of {\bf type} $Q$). 
This gives the following theorem 
(see [1], [2], [11, \S1] 
for the definition of the ``supertensor product'' 
of the superalgebras or modules).
\proclaim{Theorem 1.1}Let $C=A\lot B$ be the 
supertensor product of superalgebras $A$ and $B$ and let $V=U\otimes W$ be the 
supertensor product of a simple $A$-module $U$ and 
a simple $B$-module $W$. 
\roster 
\item"(a)" If $U$, $W$ are of type $M$, then $V$ is a simple $C$-module 
of type $M$. 
\item"(b)" If one of $U$ and $W$ is of type $M$ and the other 
is of type $Q$, then $V$ is a simple $C$-module of type $Q$. 
\item"(c)" If $U$ and $W$ are of type $Q$, then 
$V$ is a sum of two copies of a simple $C$-module $X$ 
of type $M${\rm:} $V=X\oplus X$. 
\endroster
Moreover, the above construction gives all simple $A\lot B$-modules. 
\endproclaim
Using the above $U$, $W$, $V$ and $X$, 
define an $A\lot B$-module $U\hotimes W$ by 
$$
U\hotimes W=\left\{
\matrix
\format\l&\quad\l\\
V &\t{\rm if $U$ or $W$ is of type $M$,}\\
X &\t{\rm if $U$ and $W$ are of type $Q$.}\\
\endmatrix\right. \tag1.2
$$

Let $\Irr A$ denote the set of all isomorphism classes 
of simple $A$-modules for any superalgebra $A$.
\proclaim{Corollary 1.2}We have a bijection 
$$
\hotimes\:
\Irr A\times \Irr B\; \ni (U,W)\;\overset \sim \to \mapsto \;U\hotimes W\in\;
\Irr {A\lot B}.
$$ 
\endproclaim
\comment
\subhead
1.1.4. Double centralizer theorem
\endsubsubhead
Now we introduce the double supercentralizer theorem for 
semisimple superalgebras. 
\proclaim{Theorem 5.4} \t{\rm [11, Th\. 2.1]} 
Let $A$ be a semisimple superalgebra. 
Let $V$ be an $A$-module with the associated representation
 $\rho\:A\to\E(V)$, and 
$$
V=V_1\oplus V_2\oplus\cdots\oplus V_s\tag5.4
$$
be the decomposition of $V$ into $A$-homogeneous components. 
Put $B=\E^{\cdot}_A(V)$. We define a representation 
$\T{\rho}\:A\lot B\to\E(V)$ of $A\lot B$ by 
$\T{\rho}(a\lot b)=\rho(a)\circ b$ for $a\in A$ and $b\in B$. 
Then each $V_i$, $1\leq i\leq s$, is a 
simple $A\lot B$-submodule of $V$ of type $M$. 

Let $U_i$ be a simple $A$-module contained in $V_i$.
 If $U_i$ is of type $M$ (resp\. of type $Q$), 
then there exists a simple $B$-module $W_i$ 
of type $M$ (resp\. $Q$) such that 
$$
V_i\cong_{A\lot B}U_i\hotimes W_i.
$$ 
We have $U_i\ncong_A U_j$, $W_i\ncong_B W_j$ for all $1\leq i\neq j\leq s$. 
This gives a $1-1$ correspondence between the 
isomorphism classes of simple $A$-modules 
and simple $B$-modules appearing in $V$. 
Furthermore, we have $\E^{\cdot}_B(V)=\rho(A)$. 
\endproclaim

The above theorem gives a supercommuting action of $A$ and $B$ on $V$, 
which decomposes into a multiplicity-free sum of simple 
$A\lot B$-modules, where type $M$ (resp\. of type $Q$) simple $A$-modules are 
paired with type $M$ (resp\. of type $Q$) simple $B$-modules in a bijective 
manner. 
We also encounter a similar but slightly different situation, 
in which type $M$ (resp\. of type $Q$) 
simple $A$-modules are paired with type $Q$ (resp\. of type $M$) simple 
$B$-modules. 
Here we formulate this as the following corollary. 

\proclaim{Corollary 5.5} \t{\rm[11, Cor\. 2.2]} Let $A$, $V$ and 
$V=V_1\oplus V_2\oplus\cdots\oplus V_s$ be as in Theorem 5.4. 
Assume that $\E^{1}_A(V)$ contains an invertible element $x$. 
Let $C$ denote the subsuperalgebra of $\E^{\cdot}_A(V)$ 
generated by $x$, and put $A'=A\lot C$ and $B=\E^{\cdot}_{A'}(V)$. 
Restricting the homomorphism 
$\rho\:A'\lot B\to \E(V)$ of Theorem 5.4, we can regard $V$ as an 
$A\lot B$-module. 
Then each $V_i$, $1\leq i\leq s$, is a 
simple $A\lot B$-submodule of $V$ of type $Q$. 

Let $T_i$ be a simple $A$-module contained in $V_i$.
 If $T_i$ is of type $M$ (resp\. of type $Q$), 
then there exists a simple $B$-module $W_i$ of type $Q$ 
(resp\. of type $M$) such that 
$$
V_i\cong_{A\lot B}T_i\hotimes W_i.
$$ 
We have $T_i\ncong_A T_j$, $W_i\ncong_B W_j$ for all $1\leq i\neq j\leq s$. 
This also gives a $1-1$ correspondence between the 
isomorphism classes of simple $A$-modules and simple $B$-modules appearing 
in $V$. 
Furthermore, we have 
$\E^{\cdot}_A(V)=B\lot C$ and $\E^{\cdot}_B(V)=\rho(A)\lot C$. 
\endproclaim
\endcomment
\head
\S2. Simple modules for $\B'_k$
\endhead
\comment
We will translate non-graded simple $\B'_k$-modules, 
due to J\. R\. Stembridge, 
into the language of ${\Bbb Z}_2$-graded representations. 
\endcomment
The simple $\A_k$-modules are parametrized by $DP_k$ 
(cf\. [2], [7], [9]). 
For $\lam\in DP_k$, let $V_{\lam}$ denote a simple $\A_k$-module 
indexed by $\lam$. 
Then $V_{\lam}$ is of type $M$ (resp\. of type $Q$) 
if $\lam\in DP_k^{+}$ (resp\. $\lam\in DP_k^{-}$). 
We construct a $\B'_k$-module $V_{\lam,\mu}$ 
for $(\lam,\mu)\in (DP^2)_k$ as follows. 
Define a surjective homomorphism of superalgebras 
$\pi_k\:\B'_k\to\A_k$ (resp\. $\pi_k'\:\B'_k\to\A_k$) by 
$\pi_k(\tau')=1$, $\pi_k|_{\A_k}=\id_{\A_k}$ 
(resp\.$\pi_k'(\tau')=-1$, $\pi_k'|_{\A_k}=\id_{\A_k}$). 
The simple $\A_{k'}$ (resp\. $\A_{k-k'}$)-module 
$V_{\lam}$ (resp\. $V_{\mu}$) can be lifted to a 
$\B'_{k'}$ (resp\. $\B'_{k-k'}$)-module 
via $\pi_{k'}$ (resp\. $\pi'_{k-k'}$), where $k'=|\lam|$. 
This (simple) $\B'_{k'}$ (resp\. $\B'_{k-k'}$)-module is denoted by 
$V_{\lam,\phi}$ (resp\. $V_{\phi,\mu}$). 
Let $V_{\lam,\mu}$ denote the $\B'_k$-module 
induced from the $\B'_{k'}\lot\B'_{k-k'}$-module 
$V_{\lam,\phi}\hotimes V_{\phi,\mu}$, namely 
$$
V_{\lam,\mu}=\B'_k\otimes_{\B'_{k'}\lot\B'_{k-k'}}
(V_{\lam,\phi}\hotimes V_{\phi,\mu})
$$ 
(see the definition of $\hotimes$ in \thetag{1.2}), where 
$\B'_{k'}\lot\B'_{k-k'}$ is embedded into $\B'_k$ via 
$$
\align
&\tau'\lot 1\mapsto \tau',\quad
\gamma_i\lot1\mapsto \gamma_i \quad(1\leq i\leq k'-1),\\
&1\lot\tau'\mapsto \tau'_{k'+1},\quad 
1\lot\gamma_j\mapsto\gamma_{k'+j}\quad(1\leq j\leq k-k'-1)
\endalign
$$
where $\tau'_i=\gamma_{i-1}\gamma_{i-2}\cdots\gamma_1
\tau'\gamma_1\cdots\gamma_{i-2}\gamma_{i-1}$, $1\leq i\leq k$. 

\proclaim{Theorem 2.1}{\rm\!\!\!(cf\. [10], Th\. 7.1)} 
$\{V_{\lam,\mu}\;|\;(\lam,\mu)\in (DP^2)_k\}$ is 
a complete set of the isomorphism classes of simple 
$\B'_k$-modules. $V_{\lam,\mu}$ is of type $M$ {\rm(}resp\. of type $Q${\rm)} if 
$(\lam,\mu)\in {(DP^2)}_k^{+}$ 
{\rm(}resp\. $(\lam,\mu)\in {(DP^2)}_k^{-}${\rm)} 
\endproclaim
The proof is analogous to the little group method, and is omitted. It can also 
be shown that this parametrization coincides with that by Stembridge 
in [10, Theorem 7.1] modulo the usual difference between ${\Bbb Z}_2$-graded 
and non-graded modules. 

If $(\lam,\mu)\in\!{(DP^2)}_k^{-}$, then 
fix a non-zero homogeneous element $x_{\lam,\mu}$ 
of $\E_{\B'_k}(V_{\lam,\mu})$
$\cong Q(1)$ of degree $1$. 
\head
\S3. The algebras $\B_k$ and $\C_k\lot\B'_k$
\endhead 		
For any $k\geq1$, let $\B_k$ denote the associative algebra generated by 
$\tau$ and the $\sig_i$, $1\leq i\leq k-1$, with relations 
$$
\align
&{\tau}^2=\sig_i^2=1\quad(1\leq i\leq k-1),
\quad(\sig_i\sig_{i+1})^3=1\quad(1\leq i\leq k-2),
\tag3.1\\
&(\sig_i\sig_j)^2=1\quad(|i-j|\geq2),\quad
(\tau\sig_i)^2=1\quad(2\leq i\leq k-1),\\
&(\tau\sig_1)^4=-1.\\
\endalign
$$ 
We regard $\B_k$ as a superalgebra by giving 
the generator $\tau'$ (resp\. the generator $\sig_i$, 
$1\leq i\leq k-1$) degree $1$ (resp\. degree $0$). 
The subgroup of $(\B_k)^{\times}$ 
generated by $\sig_i$, $1\leq i\leq k-1$, 
is isomorphic to the symmetric group of degree $k$ and it is denoted by 
$\fS_k$. 

Let $\C_k$ denote the $2^k$-dimensional Clifford algebra, namely 
$\C_k$ is generated by $\xi_1,\dots,\xi_k$ with relations 
$$
\xi_i^2=1,\quad \xi_i\xi_j=-\xi_j\xi_i\quad(i\neq j)\;\;.\tag3.2
$$
We regard $\C_k$ as a superalgebra 
by giving the generator $\xi_i$, $1\leq i\leq k$, degree $1$. 
$\C_k$ is a simple superalgebra. 
Let $X_k$ be a unique simple $\C_k$-module. 
If $k$ is even (resp\. odd), then $X_k$ is of type $M$ 
(resp\. of type $Q$). If $k$ is odd, then fix a non-zero element $z_k$ of 
$\E^{1}_{\C_k}(X_k)$. 

Define a linear map 
$\vth\:\B_k\to\C_k\lot\A_k$ by 
$$
\align
&\vth(\tau_i)\mapsto \xi_i\otimes1\quad\quad(1\leq i\leq k),\tag3.3\\
&\vth(\sig_j)\mapsto
\frac{1}{\sqrt{2}}(\xi_j-\xi_{j+1})\otimes\gam_j\;\;(1\leq j\leq k-1)
\endalign
$$
where $\tau_i=\sig_{i-1}\cdots\sig_1\tau\sig_1\dots\sig_{i-1}$. 
Then $\vth$ is an isomorphism of algebras (cf\. [11, Th\. 3.2]). 
For $\lam\in DP_k$, define a $\B_k$-module $W_{\lam}$ 
by $W_{\lam}=X_k\hotimes V_{\lam}$. 
By Corollary 1.2, $\{W_{\lam}\;|\;\lam\in DP_k\}$ is a complete set 
of isomorphism classes of simple $\B_k$-modules. 

Let $\hat\B_k$ denote the supertensor product (cf\. [1], [2], [11, \S1]) 
of the algebras $\C_k$ and $\B'_k$, namely $\hat\B_k=\C_k\lot\B'_k$. 
Since $\B_k\cong\C_k\lot\A_k$, 
$\B_k$ can be regarded as a subalgebra of $\hat\B_k$. 
For $(\lam,\mu)\in{(DP^2)}_k$, put 
$W_{\lam,\mu}=X_k\hotimes V_{\lam,\mu}$. 
By Theorem 1.1 and \thetag{1.2}, 
$W_{\lam,\mu}$ is of type $M$ (resp\. of type $Q$) 
if $l(\lam)+l(\mu)$ is even (resp\. odd). 
By Corollary 1.2, $\{W_{\lam,\mu}\;|\;(\lam,\mu)\in 
{(DP^2)}_k\}$ is a complete set 
of isomorphism classes of simple $\hat\B_k$-modules. 

\head
\S4. A duality of $\B'_k$ and $\fq(n_0)\oplus\fq(n_1)$
\endhead 
Let $\fq(n)$ denote the Lie subsuperalgebra 
of $\fgl(n,n)$ (denoted by $l(n,n)$ in [5]) consisting 
of the matrices of the form 
$\pmatrix 
A & B\\
B & A
\endpmatrix$. The Jacobi product 
$[\;\,,\;\,]\:\fq(n)\times\fq(n)\to\fq(n)$ is defined by 
$[X,Y]=XY-(-1)^{\overline{X}\cdot\overline{Y}}YX$ 
for homogeneous elements $X$, $Y\in \fq(n)$, 
where the symbol $\bar{\quad}$ expresses the degree 
of a homogeneous element. This Lie superalgebra is called 
the queer Lie superalgebra. 
Let $\U_n=\U(\fq(n))$ denote the universal enveloping algebra 
of $\fq(n)$, which can be regarded as a superalgebra. 
Let $W$ denote the $k$-fold supertensor product of 
the $2n$-dimensional natural representation $V=V_0\oplus V_1$, 
$\vdim{V}=(n,n)$, namely $W=V^{\otimes k}$, 
where $\vdim V$ denotes the pair $(\dim V_0,\dim V_1)$. 
We define a representation $\Theta\:\U_n\to\E(W)$ by 
$$
\Theta(X)
(v_1\otimes\cdots\otimes v_k)
=\sum_{j=1}^{k}(-1)^{\overline{X}
\cdot(\overline{v_1}+\cdots+\overline{v_{j-1}})}
v_1\otimes\cdots\otimes\overset j \to{\overset\smile\to{Xv_j}}
\otimes\cdots\otimes v_k
$$
for all homogeneous elements $X\in\fq(n)$ 
and $v_i\in V$ $(1\leq i\leq k)$. 
Note that $\U_n$ is an infinite dimensional superalgebra. 
However, for a fixed number $k$, 
$\U_n$ acts on $W$ through its finite dimensional image in $\E(W)$. 
Therefore we can use the results in \S1.1.3 
on finite dimensional superalgebras 
and their finite dimensional modules. 

Let $n_0$ and $n_1$ be two positive integers such that $n_0+n_1=n$. 
The Lie superalgebra $\fq(n_0)\oplus\fq(n_1)$ can be embedded into 
$\fq(n)$ via 
$$
\fq(n_0)\oplus\fq(n_1)\ni
\left(
\pmatrix
A&B\\
B&A
\endpmatrix,
\pmatrix
C&D\\
D&C
\endpmatrix
\right)\mapsto
\pmatrix
A&0&B&0\\
0&C&0&D\\
A&0&B&0\\
0&C&0&D\\
\endpmatrix
\in\fq(n).\tag4.1
$$

The universal enveloping algebra of 
$\fq(n_0)\oplus\fq(n_1)$ is isomorphic to 
$\U_{n_0}\lot\U_{n_1}$ which can be embedded into $\U_n$ 
as a subalgebra generated by 
the elements of $\fq(n_0)\oplus\fq(n_1)$. 

Now we define a representation 
$\Psi\:\hat\B_k\to\E(W)$, which depends on $n_0$ and $n_1$, by
$$
\align
&\Psi(\xi_i\otimes1)(v_1\otimes\cdots\otimes v_k)
=(-1)^{\overline{v_1}+\cdots+\overline{v_{i-1}}}
v_1\otimes\cdots\otimes Pv_i\otimes\cdots\otimes v_k\tag4.2\\
&\hskip8cm (1\leq i\leq k),\\
&\Psi(1\otimes\tau')(v_1\otimes\cdots\otimes v_k)
=(Qv_1)\otimes v_2\otimes\cdots\otimes v_k,\\
&\Psi(1\otimes\gamma_j)(v_1\otimes\cdots\otimes v_k)\\
&=
\frac
{(-1)^{\overline{v_1}+\cdots+\overline{v_{j-1}}}}
{\sqrt{2}}
v_1\otimes\cdots\otimes(Pv_j)\otimes v_{j+1}
\otimes\cdots\otimes v_k\\
&\hskip2cm
-\frac{(-1)^{\overline{v_1}+\cdots+\overline{v_{j-1}}+\overline{v_j}}}
{\sqrt{2}}
v_1\otimes\cdots\otimes v_i\otimes(Pv_{j+1})
\otimes\cdots\otimes v_k\\
&\hskip8cm (1\leq j\leq k-1)\\
\endalign
$$
for all homogeneous elements $v_j\in V$, $1\leq j\leq k$, 
where 
$$
\align
&P=\pmatrix
0&-\sqrt{-1}I_n\\
\sqrt{-1}I_n&0
\endpmatrix\in M(n,n)_1,\\
&Q=\pmatrix
I_{n_0}&0       &0      &0       \\
0      &-I_{n_1}&0      &0       \\
0      &0       &I_{n_0}&0       \\
0      &0       &0      &-I_{n_1} 
\endpmatrix\in Q(n)_0.
\endalign
$$
Note that, by the isomorphism $\vth\:\B_k\cong\C_k\lot\A_k\subset\hat{\B}_k$, 
$W$ can be regarded 
as a $\B_k$-module and this $\B_k$-module was investigated by 
Sergeev in [8] (cf. Theorem A). Then observing the actions 
of $\vth(\tau)$, $\vth(\sig_i)\in\B'_k$, $1\leq i\leq k-1$, on $W$, we have 
$$
\align
&\Psi\left(\vth(\tau)\right)(v_1\otimes\cdots\otimes v_k)
=(Pv_1)\otimes\cdots\otimes v_k,\tag4.3,\\
&\Psi\left(\vth(\sig_i)\right)(v_1\otimes\cdots\otimes v_k)
=(-1)^{\overline{v_i}\cdot\overline{v_{i+1}}}v_1\otimes\cdots\otimes 
v_{i+1}\otimes v_i\otimes\cdots\otimes v_k\\
\endalign
$$
for all homogeneous elements $v_j\in V$, $1\leq j\leq k-1$.

Let $W'$ be a $\U_{n_0}\lot\U_{n_1}$-submodule of $W$. 
Since $\left(\fq(n_0)\oplus\fq(n_1)\right)_0
\cong
\fgl(n_0,{\Bbb C})\oplus\fgl(n_1,{\Bbb C})$ 
as a Lie algebra, 
and $V$ is a sum of two copies 
($V_0$ and $V_1$) of the natural representation of 
$\fgl(n_0,{\Bbb C})\oplus\fgl(n_1,{\Bbb C})$, 
this embeds $W'|_{\left(\fq(n_0)\oplus\fq(n_1)\right)_0}$ 
into a sum of tensor powers of the natural 
representation, so that this representation of 
$\fgl(n_0,{\Bbb C})\oplus\fgl(n_1,{\Bbb C})$ can 
be integrated to a polynomial representation $\theta_{W'}$ of 
$GL(n_0,{\Bbb C})\times GL(n_1,{\Bbb C})$. 
Let $\Ch[W']$ denote the character of $\theta_{W'}$, 
namely 
$$
\align
&\Ch[W'](x_1,x_2,\dots,x_{n_0},y_1\cdots,y_{n_1})\\
&=\tr\theta_{W'}(\diag(x_1,x_2,\dots,x_{n_0}),\diag(y_1,y_2,\dots,y_{n_1})).
\endalign
$$ 

The following theorem determines the supercentralizer 
of $\Psi(\hat\B_k)$ in $\E(W)$ and 
describes the characters of simple $\U_{n_0}\lot\U_{n_1}$-modules appearing 
in $W$. 
\proclaim{Theorem 4.1} 
\t{\rm(1)} The two superalgebras 
$\Psi(\hat\B_k)$ and $\U_{n_0}\lot\U_{n_1}$ 
act on $W$ as the mutual supercentralizers of each other{\rm:} 
$$
\E^{\cdot}_{\Theta(\U_{n_0}\lot\U_{n_1})}(W)=\Psi(\hat\B_k),
\quad \E^{\cdot}_{\Psi(\hat\B_k)}(W)=\Theta(\U_{n_0}\lot\U_{n_1}).\tag4.4 
$$

\t{\rm(2)} The simple $\hat\B_k$-module 
$W_{\lam,\mu}$ $((\lam,\mu)\in{(DP^2)}_k)$ 
occurs in $W$ 
if and only if $l(\lam)\leq n_0$ and $l(\mu)\leq n_1$. 
Moreover we have 
$$
W\cong_{\hat\B_k\lot(\U_{n_0}\lot\U_{n_1})}
\bigoplus_{
\smallmatrix
(\lam,\mu)\in{(DP^2)}_k\\
l(\lam)\leq n_0,l(\mu)\leq n_1
\endsmallmatrix}
W_{\lam,\mu}\hotimes U_{\lam,\mu}\tag4.5
$$
where $U_{\lam,\mu}$ denotes 
a simple $\U_{n_0}\lot\U_{n_1}$-module. 

\t{\rm(3)} We have $U_{\lam,\mu}\cong_{\U_{n_0}\lot\U_{n_1}}
U_{\lam}\hotimes U_{\mu}$, where $U_{\lam}$ 
{\rm(}resp\. $U_{\mu}${\rm)} 
denotes the simple $\U_{n_0}$ 
{\rm(}resp\. $\U_{n_0}${\rm)}-module 
corresponding to the simple $\B_{|\lam|}$ 
{\rm(}resp\. $\B_{|\mu|}${\rm)}-module 
$W_{\lam}$ 
{\rm(}resp\. $W_{\mu}${\rm)} 
in Sergeev's duality 
\comment
the $|\lam|$-fold (resp\. $|\mu|$-fold) tensor product of the 
$2n_0$ (resp\. $2n_1$)-dimensional 
natural representation of $\fq(n_0)$
(resp\. $\fq(n_1)$) 
\endcomment 
{\rm(}cf\. Theorem A{\rm)}. 

\t{\rm(4)} The character values of $\Ch[U_{\lam,\mu}]$ are given as follows
{\rm:} 
$$
\align
&\Ch[U_{\lam,\mu}](x_1,x_2,\dots,x_{n_0},y_1,y_2,\dots,y_{n_1})\tag4.6\\
&\quad\quad=(\sqrt{2})^{d(\lam,\mu)-l(\lam)-l(\mu)}
Q_{\lam}(x_1,x_2,\dots,x_{n_0})Q_{\mu}(y_1,y_2,\dots,y_{n_1})
\endalign
$$
where $d\:{(DP^2)}_k\to{\Bbb Z}_2$ denotes a map defined by 
$d(\lam,\mu)=0$ {\rm(}resp\. $d(\lam,\mu)=1${\rm)} 
if $l(\lam)+l(\mu)$ is even {\rm(}resp\. $l(\lam)+l(\mu)$ is odd{\rm)}. 

\endproclaim 
\demo{Proof} First we will show the second equality of \thetag{4.4}. 
\comment
$$
\E^{\cdot}_{\Psi(\hat\B_k)}(W)=\Theta(\U_{n_0}\lot\U_{n_1}). 
$$
\endcomment
Then the first equality also follows 
from the double supercentralizer theorem (abbreviated as DSCT) 
for semisimple superalgebras (cf\. [11, Th\. 2.1]). 

By Theorem A, (1), we have 
$\E^{\cdot}_{\Psi(\vth(\B_k))}(W)\supset\Theta(\U_{n_0}\lot\U_{n_1})$, 
since $\Theta(\U_{n_0}\lot\U_{n_1})$ is a subsuperalgebra of $\Theta(\U_n)$. 
Hence $\Theta(X\otimes Y)$ commutes with $\Psi(\vth(\B_k))$ 
for any $X\in \fq(n_0)$, $Y\in\fq(n_1)$. 
By direct calculations, it can be shown that 
$\Theta(X\otimes Y)$ and $\Psi(1\otimes \tau')$ also commute. 
Since $\hat\B_k$ is generated as an algebra 
by the $\vth(\tau_i)$, $1\leq i\leq k$, 
the $\vth(\sig_j)$, $1\leq j\leq k-1$, and $1\otimes\tau'$, 
we have $\E^{\cdot}_{\Psi(\hat\B_k)}(W)\supset\Theta(\U_{n_0}\lot\U_{n_1})$. 
We need only to show that 
$$
\E^{\cdot}_{\Psi(\hat\B_k)}(W)
\subset\Theta(\U_{n_0}\lot\U_{n_1}).\tag4.7
$$
We have $\E^{\cdot}_{\Psi(\hat\B_k)}(W)\subset
\E^{\cdot}_{\Psi(\vth(\B_k))}(W)=\Theta(\U_n)$ by Theorem A, (1). 
It can be easily checked that 
$\Theta(\U_n)\subset
Q(n)\otimes\cdots\otimes Q(n)$, where $Q(n)$ denotes the underlying 
vectorspace of $\fq(n)$ (or the superalgebra it forms), so that we have 
$\E^{\cdot}_{\Psi(\hat\B_k)}(W)\subset Q(n)\otimes\cdots\otimes Q(n)$. 
We identify 
$\E(W)$ with $\oversetbrace{k}
\to{\E(V)\lot\cdots\lot\E(V)}$ by defining the action of 
$X_1\otimes X_2\otimes\cdots\otimes X_k
\in\E(V)^{\lot k}$ on $W$ by 
$$
\align
&(X_1\otimes X_2\otimes\cdots\otimes X_k)
(v_1\otimes v_2\otimes\cdots\otimes v_k)\\
&\quad=(-1)^{\overline{X_2}\cdot\overline{v_1}+
\overline{X_3}\cdot(\overline{v_1}+\overline{v_2})+\cdots+
\overline{X_k}\cdot(\overline{v_1}+\cdots+\overline{v_{k-1}})}
X_1v_1\otimes X_2v_2\otimes\cdots\otimes X_kv_k
\endalign
$$
for all homogeneous elements $X_j\in\E(V)$ and $v_j\in V$, $1\leq j\leq k$. 
Define a representation $\theta\:{\Bbb C}\fS_k\to\E((\E(W))$ 
of ${\Bbb C}\fS_k$ by 
$$
\align
&\theta(\sig_i)
(X_1\otimes\cdots\otimes X_i\otimes X_{i+1}\otimes\cdots\otimes X_k)\\
&\quad=(-1)^{\overline{X_i}\cdot\overline{X_{i+1}}}
(X_1\otimes\cdots\otimes X_{i+1}\otimes X_i\otimes\cdots\otimes X_k)\\
\endalign
$$
for $1\leq i\leq k-1$ and 
homogeneous elements $X_j$, 
$1\leq j\leq k$, of $\E(V)$. 
Moreover, define elements $T_i$, $1\leq i\leq k$, of $\E(\E(W))$ by 
$$
T_i(X_1\otimes\cdots\otimes X_k)
=X_1\otimes\cdots\otimes QX_iQ\otimes\cdots\otimes X_k
$$
for all $X_j\in\E(V)$, $1\leq j\leq k$. Furthermore, put 
$$
\align
&S=\frac{1}{n!}\sum_{w\in\fS_k}\theta(w),\\
&T=\prod_{i=1}^k\left(\frac{1}{2}(\Id_{\E(W)}+T_i)\right).
\endalign
$$ 
Note that, since $T_i\in\E^{0}(\E(W))$ for all $i$, 
the factors in the definition of $T$ commute. 
If $f\in\E^{\cdot}_{\Psi(\hat\B_k)}(W)$, then 
it follows that $S(f)=f$ and 
$\dfrac{1}{2}(\Id_{\E(W)}+T_i)(f)$
$=f$, $1\leq i\leq k$, 
since 
$\theta(\sig)(f)=\Psi(\vth(\sig))\circ f\circ \Psi(\vth(\sig))^{-1}$ and 
$T_i(f)=\Psi(1\otimes\tau'_i)\circ f\circ \Psi(1\otimes\tau'_i)$. 
Therefore, any element $f$ of 
$\E^{\cdot}_{\Psi(\hat\B_k)}(W)$ 
can be expressed as a linear 
combination of elements of the form 
$$
ST(X_1\otimes\cdots\otimes X_k) 
$$
with $X_j\in Q(n)$, $1\leq j\leq k$. 
Since 
$$
T(X_1\otimes\cdots\otimes X_k)=
\left(\frac{1}{2}\right)^k
(X_1+QX_1Q)\otimes\cdots\otimes(X_k+QX_kQ)
$$    
and $X+QXQ$ belongs to $ Q(n_0)\oplus Q(n_1)$ 
for any $X\in Q(n)$, we have 
$$
T\left( Q(n)\otimes\cdots\otimes Q(n)\right)
\subset
( Q(n_0)\oplus Q(n_1))\otimes\cdots\otimes( Q(n_0)\oplus Q(n_1)).
$$
Hence it follows that 
$$
\E^{\cdot}_{\Psi(\hat\B_k)}(W)\subset 
S\bigl(
( Q(n_0)\oplus Q(n_1))\otimes\cdots\otimes( Q(n_0)\oplus Q(n_1))
\bigr).
$$
By induction on $k$, it can be shown that 
$$
S\bigl(\oversetbrace{k}
\to{( Q(n_0)\oplus Q(n_1))\otimes\cdots\otimes( Q(n_0)\oplus Q(n_1))}
\bigr)
$$ 
is generated as an algebra by elements of the form 
$S(X\otimes1\otimes\cdots\otimes1)=\Theta(X)$ 
with $X\in\fq(n_0)\oplus\fq(n_1)$. 
Therefore \thetag{4.7} follows. 

Next we will show (2) and (3) simultaneously. 
Since $V$ is a sum of the natural representations 
$X$ and $Y$ of $\fq(n_0)$ and $\fq(n_1)$ respectively: 
$V=X\oplus Y$, where $\vdim X=(n_0,n_0)$, $\vdim Y=(n_1,n_1)$, 
$W$ can be decomposed into a sum of tensor powers of $X$ and $Y$. 
Since a $\U_{n_0}\lot\U_{n_1}$-submodule of $W$ of the form 
$\cdots\otimes X\otimes Y\otimes \cdots$ is isomorphic to that of the form 
$\cdots\otimes Y\otimes X\otimes\cdots$, we have 
$$
W\cong_{\U_{n_0}\lot\U_{n_1}}
\bigoplus_{k'=0}^{k}\left(
\oversetbrace{k'}\to{
X\otimes\cdots\otimes X}
\otimes
\oversetbrace{k-k'}\to{Y\otimes\cdots\otimes Y}\right)
^{\oplus\binom{k}{k'}}.
$$
From Theorem A, (2), we have 
$$
\align
&X^{\otimes k'}
\cong_{\B_{k'}\lot\U_{n_0}}
\bigoplus_{
\smallmatrix
\lam \in DP_k\\
l(\lam)\leq n_0
\endsmallmatrix}
W_{\lam}\hotimes U_{\lam},\tag4.8\\
&Y^{\otimes k-k'}
\cong_{\B_{k-k'}\lot\U_{n_1}}
\bigoplus_{
\smallmatrix
\mu \in DP_k\\
l(\mu)\leq n_1
\endsmallmatrix}
W_{\mu}\hotimes U_{\mu}.
\endalign
$$
Therefore, it follows that simple 
$\U_{n_0}\lot\U_{n_1}$-modules which occur in $W$ are of the form 
$U_{\lam}\hotimes U_{\mu}$, $(\lam,\mu)\in {(DP^2)}_k$, and that 
$U_{\lam}\hotimes U_{\mu}$ occurs in $W$ 
if and only if $l(\lam)\leq n_0$ and $l(\mu)\leq n_1$. 
By \thetag{4.4} and DSCT, 
$W$ can be decomposed into a sum of non-isomorphic simple 
$\hat\B_k\lot(\U_{n_0}\lot\U_{n_1})$-modules. 
In order to determine the simple $\hat\B_k$-module which is paired with 
the simple $\U_{n_0}\lot\U_{n_1}$-module $U_{\lam}\hotimes U_{\mu}$, 
we consider the $\B_{k'}\lot\B_{k-k'}$-submodule 
$\oversetbrace{k'}\to{
X\otimes\cdots\otimes X}
\otimes\oversetbrace{k-k'}\to{Y\otimes\cdots\otimes Y}$ of $W$. 
Since $\tau'_i\in\hat\B_{k'}$, $1\leq i\leq k'$ 
(resp\. $\tau'_j\in\hat\B_{k-k'}$, $1\leq j\leq k-k'$), acts on 
$X^{\otimes k'}$ (resp\. $Y^{\otimes k-k'}$) as 
$\Id_{X^{\otimes k'}}$ (resp\. $-\Id_{Y^{\otimes k-k'}}$), 
the $\B_{k'}$ (resp\. $\B_{k-k'}$)-submodule 
$W_{\lam}$ (resp\. $W_{\mu}$) of $X^{\otimes k'}$ (resp\. $Y^{\otimes k-k'}$) 
can be regarded as a $\hat\B_{k'}$ (resp. $\hat\B_{k-k'}$)-module 
and is isomorphic to $W_{\lam,\phi}$ (resp\. $W_{\phi,\mu}$). 
From \thetag{4.8}, 
a simple $\hat\B_k$-submodule of $W$ which corresponds 
to $U_{\lam}\hotimes U_{\mu}$ contains 
$W_{\lam,\phi}\otimes W_{\phi,\mu}$ 
as a $\hat\B_{k'}\lot\hat\B_{k-k'}$-submodule. 
This condition forces this simple $\hat\B_k$-module to be isomorphic 
to $W_{\lam,\mu}$. 
Consequently, the result (2) and (3) follow. 

The result (4) 
immediately follows from Theorem A, (3) and the fact that 
$$
\align
&\Ch[U\hotimes U'](x_1,\dots,x_{n_0},y_1,\dots,y_{n_1})\\
&=
\left\{
\matrix
\format\l&\,\l\\
\Ch[U](x_1,\dots,x_{n_0})\Ch[U'](y_1,\dots,y_{n_1})
&\t{if $U$ or $U'$ is of type $M$,}\\
\dfrac{1}{2}\Ch[U](x_1,\dots,x_{n_0})\Ch[U'](y_1,\dots,y_{n_1})
&\t{if $U$, $U'$ are of type $Q$.}
\endmatrix\right.
\endalign
$$
\qed\enddemo

By Theorem 4.1, (3), Theorem 1.1, \thetag{1.2} and Theorem A, 
the simple $\U_{n_0}\lot\U_{n_1}$-module $U_{\lam,\mu}$ is 
of type $M$ (resp\. of type $Q$) if $l(\lam)+l(\mu)$ is even (resp\. odd). 
If $l(\lam)+l(\mu)$ is odd, then fix a non-zero element $u_{\lam,\mu}$ 
of $\E^1_{\U_{n_0}\lot\U_{n_1}}(U_{\lam,\mu})$. 

We can rewrite \thetag{4.5} using 
the isomorphism $W_{\lam,\mu}\cong X_k\hotimes V_{\lam,\mu}$ as 
$\hat\B_k$-modules. We have 
$$
W\cong
\bigoplus_{(\lam,\mu)\in {(DP^2)}_k}
X_k\hotimes V_{\lam,\mu}\hotimes U_{\lam,\mu}.
$$
Note that, if $U$, $V$ and $W$ are simple modules for superalgebras $A$, $B$ 
and $C$ respectively, then both 
$(U\hotimes V)\hotimes W$ and $U\hotimes(V\hotimes W)$ denote the unique 
(up to isomorphism) simple $(A\lot B\lot C)$-module occurring in 
$(U\otimes V)\otimes W\cong U\otimes(V\otimes W)$, so that, up to isomorphism, 
the operation $\hotimes$ is associative. 
There are three cases where 
the $\C_k\lot\B'_k\lot(\U_{n_0}\lot\U_{n_1})$-module 
$X_k\hotimes V_{\lam,\mu}\hotimes U_{\lam,\mu}$ 
is different from the supertensor product 
$X_k\otimes V_{\lam,\mu}\otimes U_{\lam,\mu}$. 

(1) If $k$ is even and $(\lam,\mu)\in{(DP^2)}_k^{-}$, then 
$X_k$, $V_{\lam,\mu}$, $U_{\lam,\mu}$ are of type $M$, $Q$, $Q$ 
respectively. We have 
$$
X_k\hotimes V_{\lam,\mu}\hotimes U_{\lam,\mu}
=X_k\otimes(V_{\lam,\mu}\hotimes U_{\lam,\mu})
$$ 
where $V_{\lam,\mu}\hotimes U_{\lam,\mu}$ is one of the two eigenspaces of 
$x_{\lam,\mu}\otimes u_{\lam,\mu}$. 

(2) If $k$ is odd and $(\lam,\mu)\in{(DP^2)}_k^{+}$, then 
$X_k$, $V_{\lam,\mu}$, $U_{\lam,\mu}$ are of type $Q$, $M$, $Q$ 
respectively. We have 
$$
X_k\hotimes V_{\lam,\mu}\hotimes U_{\lam,\mu}
=(X_k\otimes V_{\lam,\mu})\hotimes U_{\lam,\mu}
$$ 
where $(X_k\otimes V_{\lam,\mu})\hotimes U_{\lam,\mu}$ is one 
of the two eigenspaces of $(z_k\otimes1)\otimes u_{\lam,\mu}$. 

(3) If $k$ is odd and $(\lam,\mu)\in{(DP^2)}_k^{-}$, then 
$X_k$, $V_{\lam,\mu}$, $U_{\lam,\mu}$ are of type $Q$, $Q$, $M$ 
respectively. We have 
$$
X_k\hotimes V_{\lam,\mu}\hotimes U_{\lam,\mu}
=(X_k\hotimes V_{\lam,\mu})\otimes U_{\lam,\mu}
$$ 
where $X_k\hotimes V_{\lam,\mu}$ is 
one of the two eigenspaces of $z_k\otimes x_{\lam,\mu}$. 

Put $r=\lfl k/2\rfl$ and 
$\zeta_i=\sqrt{-1}\xi_{2i-1}\xi_{2i}\in\C_k$ for 
$1\leq i\leq r$. Then the $\Psi(\zeta_i\otimes1)$, $1\leq i\leq r$, are 
commuting involutions of 
$\Psi((\C_k)_0\lot1)\subset\Psi((\hat\B_k)_0)
=\E^{0}_{\Theta(\U_{n_0}\lot\U_{n_1})}(W)$. 
For each $\vep=(\vep_1,\dots,\vep_r)\in{\Bbb Z}_2^r$, put 
$W^{\vep}
=\{w\in W\,;\,\Psi(\zeta_i\otimes1)(w)=(-1)^{\vep_i}w\quad(1\leq i\leq r)\}$. 
Then we have $W=\bigoplus_{\vep\in{\Bbb Z}_2^r}W^{\vep}$. 
Since $\zeta_i\otimes1$ commutes with $1\lot\B'_k$ for each $1\leq i\leq r$, 
$W^{\vep}$ is a $\B'_k\lot(\U_{n_0}\lot\U_{n_1})$-module. 
\proclaim{Theorem 4.2} 
For each $\vep\in{\Bbb Z}_2^r$, the submodule $W^{\vep}$ 
is decomposed as a multiplicity-free sum of simple 
$\B'_k\lot(\U_{n_0}\lot\U_{n_1})$-modules as follows{\rm:} 
$$
W^{\vep}\cong_{\B'_k\lot(\U_{n_0}\lot\U_{n_1})}
\bigoplus_{(\lam,\mu)\in {(DP^2)}_k}
V_{\lam,\mu}\hotimes U_{\lam,\mu}.
\tag4.9
$$
In the above decomposition, 
the simple $\B'_k$-modules are 
paired with the simple 
$\U_{n_0}\lot\U_{n_1}$-modules in a bijective manner. 
More precisely, we have the following results. 

\t{\rm(1)} Assume that $k$ is even. 
Then the simple $\B'_k\lot(\U_{n_0}\lot\U_{n_1})$-modules 
$V_{\lam,\mu}\hotimes U_{\lam,\mu}$ in $W^{\vep}$ are all of type $M$. 
Furthermore we have 
$$
\E^{\cdot}_{\Theta(\U_{n_0}\lot\U_{n_1})}(W^{\vep})=\Psi(\B'_k),\quad
\E^{\cdot}_{\Psi(\B'_k)}(W^{\vep})=\Theta(\U_{n_0}\lot\U_{n_1}).\tag4.10
$$

\t{\rm(2)} Assume that $k$ is odd. Then 
the simple $\B'_k\lot(\U_{n_0}\lot\U_{n_1})$-modules 
$V_{\lam,\mu}\hotimes U_{\lam,\mu}$ in $W^{\vep}$ are all of type $Q$. 
Furthermore we have 
$$
\E^{\cdot}_{\Theta(\U_{n_0}\lot\U_{n_1})}(W^{\vep})\cong\C_1\lot\Psi(\B'_k),
\quad
\E^{\cdot}_{\Psi(\B'_k)}(W^{\vep})\cong\C_1\lot\Theta(\U_{n_0}\lot\U_{n_1}).
\tag4.11 
$$
\endproclaim
\demo{Proof} For each $\vep=(\vep_1,\dots,\vep_r)\in{\Bbb Z}_2^r$, put 
$X_k^{\vep}=
\{\xi\in X_k\,;\,\zeta_i\xi=(-1)^{\vep_i}\xi
\quad(1\leq i\leq r)\}$. Then we have 
$X=\bigoplus_{\vep\in{\Bbb Z}_2^r}X_k^{\vep}$. 

(1) Assume that $k$ is even. 
Note that $X_k^{\vep}$ is one-dimensional. Let 
$\xi^{\vep}$ be a base of $X_k^{\vep}$, 
namely $X_k^{\vep}={\Bbb C}\xi^{\vep}$. 
Since the $\zeta_i$ are of degree $0$, 
$\xi^{\vep}$ is a homogeneous element of $X_k$. 
Hence we have 
$X_k^{\vep}\otimes
(V_{\lam,\mu}\hotimes U_{\lam,\mu})
\cong_{\B'_k\lot(\U_{n_0}\lot\U_{n_1})}
V_{\lam,\mu}\hotimes U_{\lam,\mu}$. 

If $(\lam,\mu)\in{(DP^2)}_k^{+}$, then we have 
$$
\align
X_k\hotimes V_{\lam,\mu}\hotimes U_{\lam,\mu}
=X_k\otimes V_{\lam,\mu}\otimes U_{\lam,\mu}
&=\bigoplus_{\vep\in{\Bbb Z}_2^r}
X_k^{\vep}\otimes V_{\lam,\mu}\otimes U_{\lam,\mu}\\
&=\bigoplus_{\vep\in{\Bbb Z}_2^r}
X_k^{\vep}\otimes(V_{\lam,\mu}\hotimes U_{\lam,\mu}).\\
\endalign
$$
If $(\lam,\mu)\in{(DP^2)}_k^{-}$, then we have 
$$
X_k\hotimes V_{\lam,\mu}\hotimes U_{\lam,\mu}
=\bigoplus_{\vep\in{\Bbb Z}_2^r}
X_k^{\vep}\otimes
(V_{\lam,\mu}\hotimes U_{\lam,\mu})
$$
since the $\zeta_i$, $1\leq i\leq r$, 
and $1\otimes x_{\lam,\mu}\otimes u_{\lam,\mu}$ commute. 
Consequently we have 
$$
\align
W^{\vep}&
\cong_{\B'_k\lot(\U_{n_0}\lot\U_{n_1})}
X_k^{\vep}\otimes\left(
\bigoplus_{(\lam,\mu)\in{(DP^2)}_k}
V_{\lam,\mu}\hotimes U_{\lam,\mu}
\right)\\
&\cong_{\B'_k\lot(\U_{n_0}\lot\U_{n_1})}
\bigoplus_{(\lam,\mu)\in{(DP^2)}_k}
V_{\lam,\mu}\hotimes U_{\lam,\mu}
\endalign
$$
Therefore \thetag{4.9} follows. 
By Theorem 1.1 and \thetag{1.2}, the simple modules 
$V_{\lam,\mu}\hotimes U_{\lam,\mu}$ appearing in the above 
decomposition are of type $M$. 

First we will show the second equality in \thetag{4.10}. 
Then the first equality follows from DSCT. 
Since $W^{\vep}$ is a $\B'_k\lot(\U_{n_0}\lot\U_{n_1})$-module, we have 
$$
\Theta(\U_{n_0}\lot\U_{n_1})
|_{W^{\vep}}\subset\E^{\cdot}_{\Psi(\B'_k)}(W^{\vep}).
$$ 
By DSCT, \thetag{4.9} (already proved for this case) 
and \thetag{4.5}, we have 
$$
\dim\E^{\cdot}_{\Psi(\B'_k)}(W^{\vep})
=\dim\E^{\cdot}_{\Psi(\hat\B_k)}(W)
$$
since both equal 
$\sum\limits_{(\lam,\mu)\in{(DP^2)}_k^{+}}(\dim U_{\lam,\mu})^2+
\sum\limits_{(\lam,\mu)\in{(DP^2)}_k^{-}}\frac{1}{2}(\dim U_{\lam,\mu})^2$. 
By Theorem 4.1, (1), we have 
$\dim\E^{\cdot}_{\Psi(\hat\B_k)}(W)=\dim\Theta(\U_{n_0}\lot\U_{n_1})$. 
Define a linear map $\pr_{\vep}\:
\Theta(\U_{n_0}\lot\U_{n_1})\to
\Theta(\U_{n_0}\lot\U_{n_1})|_{W^{\vep}}$ by 
$\pr_{\vep}(f)=f|_{W^{\vep}}$ for $f\in\Theta(\U_{n_0}\lot\U_{n_1})$. 
It is clear that $\pr_{\vep}$ is surjective. 
We claim that $\pr_{\vep}$ is injective. 
Assume that $f\in\ker\pr_{\vep}$, namely 
$f|_{W^{\vep}}=0\in\E(W^{\vep})$. 
Since $f$ and the $\xi_{2j-1}$ commute, and 
a subgroup of $(\C_k)^{\times}$ generated by 
the $\xi_{2j-1}$, $1\leq j\leq r$, transitively act 
on $\{W^{\vep'}\,;\,\vep'\in{\Bbb Z}_2^r\}$ as follows{\rm:}
$$
\xi_{2j-1}W^{(\vep_1,\dots,\vep_r)}
=W^{(\vep_1,\dots,\vep_j+1,\dots,\vep_r)}\quad(1\leq^\forall j\leq r) 
$$
it follows that $f|_{W^{\vep'}}=0$ for all $\vep'\in{\Bbb Z}_2^r$. 
Therefore $f=0$ in $\E(W)$. Hence $\pr_{\vep}$ is injective. 
Consequently we have $\dim\Theta(\U_{n_0}\lot\U_{n_1})|_{W^{\vep}}
=\dim\E^{\cdot}_{\Psi(\B'_k)}(W^{\vep})$. 
It follows that $\Theta(\U_{n_0}\lot\U_{n_1})|_{W^{\vep}}
=\E^{\cdot}_{\Psi(\B'_k)}(W^{\vep})$, as required. 

(2) Assume that $k$ is odd. Note that $X_k^{\vep}$ is $2$-dimensional. 
Then 
$X_k^{\vep}={\Bbb C}\xi^{\vep}\oplus{\Bbb C}z_k\xi^{\vep}$.

If $(\lam,\mu)\in{(DP^2)}_k^{+}$, 
then $V_{\lam,\mu}\hotimes X_k=V_{\lam,\mu}\otimes X_k$ and 
we regard the $\B'_k\lot\C_k$-module 
$V_{\lam,\mu}\otimes X_k$ 
as a $\C_k\lot\B'_k$-module via 
$\omega_{\C_k,\B'_k}$, where 
$\omega_{\C_k,\B'_k}\:\C_k\lot\B'_k\to\B'_k\lot\C_k$ denotes 
an isomorphism of superalgebras determined by 
$\omega_{\C_k,\B'_k}(a\otimes b)=(-1)^{\bar{a}\cdot\bar{b}}b\otimes a$ for 
all homogeneous elements $a\in\C_k$ and $b\in\B'_k$. An isomorphism 
$\theta\:X_k\otimes V_{\lam,\mu}
\overset{\sim}\to{\to}V_{\lam,\mu}\otimes X_k$ is defined by 
$\theta(\xi\otimes v)=(-1)^{\overline{\xi}\cdot\overline{v}}v\otimes \xi$ 
for all homogeneous elements $\xi\in X_k$ and $v\in V_{\lam,\mu}$. 
Since $\theta\circ(z_k\otimes1)=(1\otimes z_k)\circ\theta$, we have 
$$
\align
X_k\hotimes V_{\lam,\mu}\hotimes U_{\lam,\mu}
&\cong_{\hat\B_k\lot\U_n}
(V_{\lam,\mu}\otimes X_k)\hotimes U_{\lam,\mu}\\
&\cong_{\hat\B_k\lot\U_n}
V_{\lam,\mu}\otimes(X_k\hotimes U_{\lam,\mu})\\
\endalign
$$
where $X_k\hotimes U_{\lam,\mu}$ 
denotes one of the two eigenspaces of 
$z_k\otimes u_{\lam,\mu}$. 
Since the $\zeta_i$ and $z_k\otimes u_{\lam,\mu}$ commute, we have 
$$
X_k\hotimes U_{\lam,\mu}=
\bigoplus_{\vep\in{\Bbb Z}_2^r}
X_k^{\vep}\hotimes U_{\lam,\mu}
$$
where $X_k^{\vep}\hotimes U_{\lam,\mu}$ 
denotes one of the two eigenspaces of 
$z_k|_{X_k^{\vep}}\otimes u_{\lam,\mu}$. 
Since $X_k^{\vep}\hotimes U_{\lam,\mu}$ is a 
$\U_{n_0}\lot\U_{n_1}$-submodule of 
$X_k\hotimes U_{\lam,\mu}\cong_{\U_{n_0}\lot\U_{n_1}}
U_{\lam,\mu}^{\oplus 2^r}$ and 
$\dim(X_k^{\vep}\hotimes U_{\lam,\mu})=\dim U_{\lam,\mu}$, 
it follows that $X_k^{\vep}\hotimes U_{\lam,\mu}
\cong_{\U_{n_0}\lot\U_{n_1}}U_{\lam,\mu}$. 

If $(\lam,\mu)\in{(DP^2)}_k^{-}$, then we have 
$$
X_k\hotimes V_{\lam,\mu}\hotimes U_{\lam,\mu}
=(X_k\hotimes V_{\lam,\mu})\otimes U_{\lam,\mu}
=\bigoplus_{\vep\in{\Bbb Z}_2^r}
(X_k^{\vep}\hotimes V_{\lam,\mu})\otimes U_{\lam,\mu}
$$
since the $\zeta_i$, $1\leq i\leq r$, and 
$z_k\otimes x_{\lam,\mu}\otimes1$ commute, 
where $X_k^{\vep}\hotimes V_{\lam,\mu}$ 
denotes one of the two eigenspaces of 
$z_k|_{X_k^{\vep}}\otimes x_{\lam,\mu}$. 
Since $X_k^{\vep}\hotimes V_{\lam,\mu}$ is a 
$\B'_k$-submodule of 
$X_k\hotimes V_{\lam,\mu}\cong_{\B'_k}
V_{\lam,\mu}^{\oplus 2^r}$ and 
$\dim(X_k^{\vep}\hotimes V_{\lam,\mu})=\dim V_{\lam,\mu}$, 
it follows that 
$X_k^{\vep}\hotimes V_{\lam,\mu}\cong_{\B'_k} V_{\lam,\mu}$. 

Consequently we have 
$$
W^{\vep}\cong\bigoplus_{(\lam,\mu)\in{(DP^2)}_k}
V_{\lam,\mu}\otimes U_{\lam,\mu}.
$$
By Theorem 1.1 and \thetag{1.2}, the simple modules 
$V_{\lam,\mu}\otimes U_{\lam,\mu}$ appearing in the above decomposition 
are of type $Q$ and we have 
$V_{\lam,\mu}\otimes U_{\lam,\mu}=V_{\lam,\mu}\hotimes U_{\lam,\mu}$. 
Therefore, \thetag{4.9} and the former statement of (2) follows. 

The supercentralizer 
$\E^{\cdot}_{\Psi(\B'_k)}(W^{\vep})$ contains 
an invertible element $\Psi(\xi_k)\hskip-.08cm\in\hskip-.1cm
\Psi(\C_k)$. The subsuperalgebra of 
$\E^{\cdot}_{\Psi(\B'_k)}(W^{\vep})$ generated by $\Psi(\xi_k)$ is 
isomorphic to $\C_1$. By the arguments similar to the proof of \thetag{4.10}, 
the result \thetag{4.11} follows from DSCT (cf\. [11, Cor\. 2.1]). 
\qed
\enddemo

Let us mention a relation between the branching rule of the $\fq(n)$-modules 
to $\fq(n_0)\oplus\fq(n_1)$ and 
that of the $\hat\B_k$-modules to $\B_k$ 
(or that of the $\B'_k$-modules to $\A_k$). 

If an $A$-module $V$ restricts to an $B$-module, we write 
$V\downarrow^{A}_{B}$ for this $B$-module, for a superalgebra 
$A$ and a subsuperalgebra $B$ of $A$. Moreover, we write $[V:U]_{A}$ 
(or simply write $[V:U]$) for the multiplicity of a simple $A$-module 
$U$ in an $A$-module $V$. 
\proclaim{Corollary 4.3} Put 
$$
\align
&m_{\mu,\nu}^{\lam}=
[U_{\lam}\downarrow^{\U_n}_{\U_{n_0}\lot\U_{n_1}}
:U_{\mu,\nu}],\tag4.12\\
&{m'}_{\mu,\nu}^{\lam}=
[W_{\mu,\nu}\downarrow^{\hat\B_k}_{\B_k}
:W_{\lam}]
\quad
\left(\t{\rm resp\.} 
[V_{\mu,\nu}
\downarrow^{\B'_k}_{\A_k}
:V_{\lam}
]\right).
\tag4.13
\endalign
$$
Then we have 
$$
\align
&{m'}_{\mu,\nu}^{\lam}=\left\{
\matrix
\format\l&\quad\l\\
\dfrac{1}{2}m_{\mu,\nu}^{\lam}
&\t{if $W_{\mu,\nu}$ (resp\. $V_{\mu,\nu}$) is of type $M$}\\ 
&\t{\quad  and $W_{\lam}$ (resp\. $V_{\lam}$) is of type $Q$,}\\
2m_{\mu,\nu}^{\lam}
&\t{if $W_{\mu,\nu}$ (resp\. $V_{\mu,\nu}$) is of type $Q$}\\ 
&\t{\quad  and $W_{\lam}$ (resp\. $V_{\lam}$) is of type $M$,}\\
m_{\mu,\nu}^{\lam}
&\t{otherwise.}\\
\endmatrix\right.\tag4.14\\
\endalign
$$
\endproclaim
\demo{Proof} Put 
$$
W'=W_{\lam}\hotimes U_{\lam,\mu},\quad 
W_1
=W\downarrow^{\B_k\lot\U_n}_{\B_k\lot(\U_{n_0}\lot\U_{n_1})},\quad
W_2
=W\downarrow^{\hat\B_k\lot(\U_{n_0}\lot\U_{n_1})}
_{\B_k\lot(\U_{n_0}\lot\U_{n_1})}.
$$ 
Since $W_1\cong W_2$, we have $[W_1:W']=[W_2:W']$. 
Moreover, put 
$$
W'_1=
(W_{\lam}\hotimes U_{\lam})\downarrow
^{\B_k\lot\U_n}_{\B_k\lot(\U_{n_0}\lot\U_{n_1})},\,
W'_2=(W_{\mu,\nu}\hotimes U_{\mu,\nu})
\downarrow^{\hat\B_k\lot(\U_{n_0}\lot\U_{n_1})}
_{\B_k\lot(\U_{n_0}\lot\U_{n_1})}.
$$ 
From \thetag{4.5} and \thetag{A.2}, we have 
$[W_1:W']=[W'_1:W']$ and $[W_2:W']=[W'_2:W']$. 
Using \thetag{4.12} and \thetag{4.13}, we have 
$$
\align
&
(W_{\lam}\otimes U_{\lam})\downarrow
^{\B_k\lot\U_n}_{\B_k\lot(\U_{n_0}\lot\U_{n_1})}
\cong\bigoplus_{(\mu,\nu)\in{(DP^2)}_k}
\left(W_{\lam}\otimes U_{\mu,\nu}
\right)^{\oplus m_{\mu,\nu}^{\lam}},\\
&(W_{\mu,\nu}\otimes U_{\mu,\nu})
\downarrow^{\hat\B_k\lot(\U_{n_0}\lot\U_{n_1})}
_{\B_k\lot(\U_{n_0}\lot\U_{n_1})}
\cong\bigoplus_{\lam\in DP_k}
\left(W_{\lam}\otimes U_{\mu,\nu}
\right)^{\oplus {m'}_{\mu,\nu}^{\lam}}.\\
\endalign
$$
By \thetag{1.2}, the above supertensor products 
$(W_{\lam}\otimes U_{\lam})\downarrow
^{\B_k\lot\U_n}_{\B_k\lot(\U_{n_0}\lot\U_{n_1})}$, 
$(W_{\mu,\nu}\otimes U_{\mu,\nu})
\downarrow^{\hat\B_k\lot(\U_{n_0}\lot\U_{n_1})}
_{\B_k\lot(\U_{n_0}\lot\U_{n_1})}$, $W_{\lam}\otimes U_{\mu,\nu}$ 
are sums of two copies of 
$W'_1$, $W'_2$, $W'$ if $l(\mu)+l(\nu)$ is odd, $l(\lam)$ is odd, 
$l(\lam)$ and $l(\mu)+l(\nu)$ are odd, respectively. 
Therefore we have 
$$
\align
&[W'_1:W']=\left\{
\matrix
\format\l&\quad\l\\
\dfrac{1}{2}m_{\mu,\nu}^{\lam}
&\t{if $l(\lam)$ is odd and $l(\mu)+l(\nu)$ is even,}\\
m_{\mu,\nu}^{\lam}&\t{ otherwise,}\\
\endmatrix\right.,\\
&[W'_2:W']=\left\{
\matrix
\format\l&\quad\l\\
\dfrac{1}{2}{m'}_{\mu,\nu}^{\lam}
&\t{ if $l(\lam)$ is even and $l(\mu)+l(\nu)$ is odd,}\\
{m'}_{\mu,\nu}^{\lam}&\t{ otherwise.}\\
\endmatrix\right.
\endalign
$$
Comparing the above two equations, the result 
\thetag{4.14} follows. 

Next, using \thetag{4.9} and \thetag{B.1},
 we consider the multiplicities of the simple 
$\A_k\lot(\U_{n_0}\lot\U_{n_1})$-module 
$V_{\lam}\hotimes U_{\lam,\mu}$ in  
$W^{\vep}\downarrow^{\A_k\lot\U_n}_{\A_k\lot(\U_{n_0}\lot\U_{n_1})}$ 
and $W^{\vep}\downarrow^{\B'_k\lot(\U_{n_0}\lot\U_{n_1})}
_{\A_k\lot(\U_{n_0}\lot\U_{n_1})}$ respectively. 
Then \thetag{4.14} similarly follows. 
\comment
Since 
$W\downarrow^{\hat\B_k\lot(\U_{n_0}\lot\U_{n_1})}
_{\B_k\lot(\U_{n_0}\lot\U_{n_1})}
\cong_{\B_k\lot(\U_{n_0}\lot\U_{n_1})}
W\downarrow^{\B_k\lot\U_n}_{\B_k\lot(\U_{n_0}\lot\U_{n_1})}$, we have 
From \thetag{4.*}, we have  
$$
\align
&(V_{\lam}\otimes U_{\lam})
\downarrow^{\A_k\lot\U_n}_{\A_k\lot(\U_{n_0}\lot\U_{n_1})}
\cong\bigoplus_{(\mu,\nu)\in{(DP^2)}_k}
\left(V_{\lam}\otimes U_{\mu,\nu}
\right)^{\oplus m_{\mu,\nu}^{\lam}}\\
&(V_{\mu,\nu}\otimes U_{\mu,\nu})
\downarrow^{\B'_k\lot(\U_{n_0}\lot\U_{n_1})}
_{\A_k\lot(\U_{n_0}\lot\U_{n_1})}
\cong\bigoplus_{\lam\in DP_k}
\left(V_{\lam}\otimes U_{\mu,\nu}
\right)^{\oplus {m''}_{\mu,\nu}^{\lam}}
\endalign
$$
Note that 
$$
\align
&\hskip2cm V_{\lam}\otimes U_{\lam}\cong\left\{
\matrix
\format\l&\quad\l\\
(V_{\lam}\hotimes U_{\lam})^{\oplus 2}
&\t{ if $k$ is even and $\lam\in DP_k^{-}$,}\\
V_{\lam}\hotimes U_{\lam}&\t{ otherwise,}
\endmatrix\right.\\
&\left(\t{resp\. }
V_{\mu,\nu}\otimes U_{\mu,\nu}
\cong\left\{
\matrix
\format\l\\
(V_{\mu,\nu}\hotimes U_{\mu,\nu})^{\oplus 2}\\
\quad\quad\quad
\t{if $k$ is even and $(\mu,\nu)\in{(DP^2)}_k^{-}$,}\\
V_{\mu,\nu}\hotimes U_{\mu,\nu}\\
\quad\quad\quad\t{otherwise,}
\endmatrix\right.\right)
\endalign
$$
as $\A_k\lot\U_n$ (resp\. $\B'_k\lot(\U_{n_0}\lot\U_{n_1})$)-modules, 
and therefore as $\A_k\lot(\U_{n_0}\lot\U_{n_1})$-modules, and 
$$
V_{\lam}\otimes U_{\mu,\nu}\cong\left\{
\matrix
\format\l&\quad\l\\
(V_{\lam}\hotimes U_{\mu,\nu})^{\oplus 2}
&\t{ if $\lam\in DP_k^{-}$, $l(\mu)+l(\nu)$ is odd,}\\
V_{\lam}\hotimes U_{\mu,\nu}
&\t{ otherwise}\\
\endmatrix\right.
$$
as $\A_k\lot(\U_{n_0}\lot\U_{n_1})$-modules. 
Therefore, multiplicities of simple $\A_k\lot(\U_{n_0}\lot\U_{n_1})$-modules 
$V_{\lam}\hotimes U_{\mu,\nu}$ in 
$(V_{\lam}\hotimes U_{\lam})
\downarrow^{\A_k\lot\U_n}_{\A_k\lot(\U_{n_0}\lot\U_{n_1})}$ (resp\. 
$V_{\mu,\nu}\hotimes U_{\mu,\nu}
\downarrow^{\B'_k\lot(\U_{n_0}\lot\U_{n_1})}
_{\A_k\lot(\U_{n_0}\lot\U_{n_1})}$)
are as follows{\rm:} 
$$
\align
&[V_{\lam}\hotimes U_{\mu,\nu}:
V_{\lam}\hotimes U_{\lam}
\downarrow^{\A_k\lot\U_n}_{\A_k\lot(\U_{n_0}\lot\U_{n_1})}
]_{\A_k\lot(\U_{n_0}\lot\U_{n_1})}\\
&\quad=\left\{
\matrix
\format\l&\quad\l\\
\dfrac{1}{2}m_{\mu,\nu}^{\lam}
&\t{if $k$ is even and $\lam\in DP_k^{-}$, $(\mu,\nu)\in{(DP^2)}_k^{+}$}\\
2m_{\mu,\nu}^{\lam}
&\t{if $k$ is odd and $\lam\in DP_k^{-}$, $(\mu,\nu)\in{(DP^2)}_k^{+}$}\\
m_{\mu,\nu}^{\lam}
&\t{ otherwise}\\
\endmatrix\right.\\
\intertext{and }
&[
V_{\lam}\hotimes U_{\mu,\nu}:
V_{\mu,\nu}\hotimes U_{\mu,\nu}
\downarrow^{\B'_k\lot(\U_{n_0}\lot\U_{n_1})}
_{\A_k\lot(\U_{n_0}\lot\U_{n_1})}
]_{\A_k\lot(\U_{n_0}\lot\U_{n_1})}
\\
&\quad=\left\{
\matrix
\format\l&\quad\l\\
\dfrac{1}{2}{m''}_{\mu,\nu}^{\lam}
&\t{ if $k$ is even and $\lam\in DP_k^{+}$ 
, $(\mu,\nu)\in{(DP^2)}_k^{-}$}\\
2{m''}_{\mu,\nu}^{\lam}
&\t{ if $k$ is odd and $\lam\in DP_k^{-}$ 
, $(\mu,\nu)\in{(DP^2)}_k^{+}$}\\
{m''}_{\mu,\nu}^{\lam}&\t{ otherwise.}\\
\endmatrix\right.
\endalign
$$
Therefore \thetag{4.*} follows. 
\endcomment
\qed\enddemo 

Let $H'_k$ be the subgroup of $(\B'_k)^{\times}$ 
generated by $-1$, $\tau'$, $\gamma_1,\dots,\gamma_{k-1}$. 
Then $H'_k$ is a double cover 
(a central extension with a ${\Bbb Z}_2$ kernel) of $H_k$. 
Let $w^{\kappa,\nu}$ denote the element of $H'_k$ 
defined by 
$$
\align
&w^{\kappa,\nu}=w_1w_2\cdots w_{l}w'_1w'_2\cdots w'_{l'}\quad
(l=l(\kappa),l'=l(\nu)),\\
&w_i=\gam_{a+1}\gam_{a+2}\cdots\gam_{a+\kappa_i-1}\quad
(a=\kappa_1+\cdots+\kappa_{i-1}),\\
&w'_i=\gam_{b+1}\gam_{b+2}\cdots\gam_{b+\nu_i-1}\tau'_{b+\nu_i}\quad
(b=|\kappa|+\nu_1+\cdots+\nu_{i-1}).\\
\endalign
$$
Note that the image of $w^{\kappa,\nu}$ in $H_k$ is 
a representative of the conjugacy class of $H_k$ indexed by $(\kappa,\nu)$. 

Define a map $\vep\:{(DP^2)}_k\to{\Bbb Z}_2$ by 
$\vep(\lam,\mu)=1$ (resp\. $\vep(\lam,\mu)=0$) 
if $(\lam,\mu)\in{(DP^2)}_k^{+}$ 
(resp\. $(\lam,\mu)\in{(DP^2)}_k^{-}$). 

We describe a formula for character values of simple $\B'_k$-modules. 
\proclaim{Corollary 4.4} We have 
$$
\align
&2^{\frac{l(\kappa)+l(\nu)}{2}}
p_{\kappa}(x,y)p_{\nu}(x,-y)\tag4.15\\
&\quad\quad=
\sum_{(\lam,\mu)\in{(DP^2)}_k}
\Ch[V_{\lam,\mu}](w^{\kappa,\nu})
2^{\frac{-l(\lam)-l(\mu)-\vep(\lam,\mu)}{2}}
Q_{\lam}(x)Q_{\mu}(y)
\endalign
$$
for all $(\kappa,\nu)\in{(OP^2)}_k$, where 
$p_{\kappa}(x,y)=p_{\kappa}(x_1,x_2,\dots,y_1,y_2,\dots)$ and 
$p_{\nu}(x,y)=p_{\nu}(x_1,x_2,\dots,-y_1,-y_2,\dots)$ and 
$\Ch[V_{\lam,\mu}]$ denotes the character of $V_{\lam,\mu}$, namely 
$\Ch[V_{\lam,\mu}](w)=\tr(w_{V_{\lam,\mu}})$ for $w\in\B'_k$ where 
$w_{V_{\lam,\mu}}$ denotes the action of $w\in\B'_k$ 
on $V_{\lam,\mu}$. 
\endproclaim
\demo{Proof} By what we noted before Theorem 4.1, any 
$\B'_k\lot(\U_{n_0}\lot\U_{n_1})$-submodule $W'$ 
of $W$ can be regarded as a $\B'_k$-module with a commuting polynomial 
representation $\theta_{W'}$ of 
$GL(n_0,{\Bbb C})\times GL(n_1,{\Bbb C})$. Here we extend our 
notation in Theorem 4.1 to let 
$\Ch[W'](x\otimes g)$ denote 
the trace $\tr(x_{W'}\circ\theta_{W'}(g))$ for $x\in\B'_k$ and 
$g\in GL(n_0,{\Bbb C})\times GL(n_1,{\Bbb C})$, 
where $x_{W'}$ denotes the action of $x\in\B'_k$ on $W'$. 

For any $\vep$, $\vep'\in{\Bbb Z}_2^r$, 
we have $W^{\vep}\cong_{\B'_k\lot(\U_{n_0}\lot\U_{n_1})}W^{\vep'}$. 
Hence, for $(\kappa,\nu)\in{(OP^2)}_k$ and 
$E=\diag(x_1,\dots,x_{n_0},y_1,\dots,y_{n_1})
\in GL(n,{\Bbb C})$, we have 
$$
\Ch[W^{\vep}]\left(w^{\kappa,\nu}\otimes E\right)
=2^{-r}
\Ch[W]\left((1\otimes w^{\kappa,\nu})\otimes E\right)\tag4.16
$$
where $1\otimes w^{\kappa,\nu}\in\C_k\lot\B'_k=\hat\B_k$. 
We calculate the right hand side using 
the embedding $\vth\:\B_k\hookrightarrow\hat\B_k$ 
(cf\. \thetag{3.3}), namely 
$1\otimes\gam_j=\vth(\frac{1}{\sqrt{2}}(\tau_j-\tau_{j+1})\sig_j)$. 
Put $k'=|\kappa|$ and $l=l(\kappa)$. Then $k-k'=|\nu|$. Moreover put 
$W'=V^{\otimes k'}$ and $W''=V^{\otimes k-k'}$. 
We have $w^{\kappa,\nu}=w^{\kappa,\phi}w^{\phi,\nu}$, where 
$w^{\kappa,\phi}\in\B'_{k'}$, $w^{\phi,\nu}\in\B'_{k-k'}$. 
Define a representations of 
$\hat\B_{k'}$ on $W'$ 
(resp\. a representation of $\hat\B_{k-k'}$ on $W''$) 
by the same manner as the representation $\Psi$ of $\hat\B_k$ in $W$. 
Then the action of $1\otimes w^{\kappa,\phi}$ 
(resp\. the action of $1\otimes w^{\phi,\nu}$) on $W$ can be expressed as 
(the action of $1\otimes w^{\kappa,\phi}$ on $W'$)
$\otimes\oversetbrace{k-k'}\to{\id\otimes\cdots\otimes\id}$ (resp\. 
$\oversetbrace{k'}\to{\id\otimes\cdots\otimes\id}$ 
(the action of $1\otimes w^{\phi,\nu}$ on $W''$)). Hence we have 
$$
\align
&\Ch[W]\left((1\otimes w^{\kappa,\nu})\otimes E\right)\\
&=\Ch[W']\left((1\otimes w^{\kappa,\phi})\otimes E\right)
\Ch[W'']\left((1\otimes w^{\phi,\nu})\otimes E\right).\tag4.17
\endalign
$$
The element $1\otimes w^{\kappa,\phi}$ of $\hat\B_{k'}$ 
is a product of $k'-l$ elements 
$1\otimes\gamma_j=\vth(\frac{1}{\sqrt{2}}(\tau_j-\tau_{j+1})\sig_j)$. 
This product can be rearranged into the following form{\rm:} 
$$
(\t{constant})
\times(\t{a product of the $\vth(\tau_{p})-\vth(\tau_{q})$})\times 
(\t{a product of the $\vth(\sig_j)$}).
$$
The product of the $\vth(\sig_j)$ equals $\vth(\sig^{\kappa,\phi})$. 
Expanding the product of $\vth(\tau_p)-\vth(\tau_q)$ 
into a sum of $2^{k'-l}$ elements, we have 
$$
1\otimes w^{\kappa,\phi}
=\left(\frac{1}{\sqrt{2}}\right)^{k'-l}\times \sum 
(\t{a product of the $\vth(\tau_{p})$})\times \vth(\sig^{\kappa,\phi})
$$
where 
$$
\align
&{\sig}^{\kappa,\phi}=g_1g_2\cdots g_{l},\\ 
&g_i=\sig_{a+1}\sig_{a+2}\dots\sig_{a+\nu_i-1},\quad 
(a=\sum_{j=1}^{i-1}\kappa_j).
\endalign
$$ 
Then all terms in the summation are conjugate to 
$\vth(\sig^{\kappa,\phi})$ in $\vth((\B_k)^{\times})$. 
Therefore we have 
$$
\align
\Ch[W']\left(
(1\otimes w^{\kappa,\phi})\otimes E\right)
&=2^{k'-l}(\sqrt{2})^{l-k'}
\Ch[W']\left(\vth(\sig^{\kappa,\phi})\otimes E\right)\\
&=(\sqrt{2})^{k'+l}
p_{\kappa}(x_1,x_2,\dots,y_1,y_2,\dots).
\endalign
$$

Put $l'=l(\nu)$. Similarly we have
$$
\align
\Ch[W'']\left(
(1\otimes w^{\phi,\nu})\otimes E\right)
&=2^{k-k'-l'}(\sqrt{2})^{l'-k+k'}\Ch[W'']
(\vth({\sig'}^{\phi,\nu})\otimes E)\\
&=(\sqrt{2})^{k-k'+l'}
p_{\nu}(x_1,x_2,\dots,-y_1,-y_2,\dots)\\
\endalign
$$
where 
$$
\align
&{\sig'}^{\phi,\nu}=g'_1g'_2\cdots g'_{l'},\\ 
&g'_i=\sig_{b+1}\sig_{b+2}\dots\sig_{b+\nu_i-1}\tau'_{b+\nu_i},
\quad(b=\sum_{j=1}^{i-1}\nu_j).
\endalign
$$
By \thetag{4.16} and \thetag{4.17}, we have 
$$
\align
&\Ch[W^{\vep}]\left((1\otimes w^{\kappa,\nu})\otimes E\right)\\
&=\left\{
\matrix
\format\l\\
(\sqrt{2})^{l+l'}
p_{\kappa}(x_1,x_2,\dots,y_1,y_2,\dots)
p_{\nu}(x_1,x_2,\dots,-y_1,-y_2,\dots)\\
\hskip8cm\t{if $k$ is even,}\\
(\sqrt{2})^{l+l'+1}
p_{\kappa}(x_1,x_2,\dots,y_1,y_2,\dots)
p_{\nu}(x_1,x_2,\dots,-y_1,-y_2,\dots)\\
\hskip8cm\t{if $k$ is odd.}\\
\endmatrix\right.
\endalign
$$
On the other hand, by \thetag{4.6} and \thetag{4.9}, 
if $k$ is even, then we have 
$$
\align
&\Ch[\bigoplus_{(\lam,\mu)\in{(DP^2)}_k}
V_{\lam,\mu}\hotimes U_{\lam,\mu}]
\left(w^{\kappa,\nu}\otimes E\right)\\
&=\sum_{(\lam,\mu)\in{(DP^2)}_k}
\Ch[V_{\lam,\mu}](w^{\kappa,\nu})\\
&\hskip2cm\times
(\sqrt{2})^{-\vep(\lam,\mu)-l(\lam)-l(\mu)}
Q_{\lam}(x_1,\dots,x_{n_0})Q_{\mu}(y_1,\dots,y_{n_1}),\\
\endalign
$$
and if $k$ is odd, then we have 
$$
\align
&\Ch[\bigoplus_{(\lam,\mu)\in{(DP^2)}_k}
V_{\lam,\mu}\hotimes U_{\lam,\mu}]
\left(w^{\kappa,\nu}\otimes E\right)\\
&=\sqrt{2}\sum_{(\lam,\mu)\in{(DP^2)}_k}
\Ch[V_{\lam,\mu}](w^{\kappa,\nu})\\
&\hskip2cm\times(\sqrt{2})^{-\vep(\lam,\mu)-l(\lam)-l(\mu)}
Q_{\lam}(x_1,\dots,x_{n_0})Q_{\mu}(y_1,\dots,y_{n_1}).\\
\endalign
$$
Since these hold for all $n_0$ and $n_1$, the result follows.
\qed\enddemo  
We review Stembridge's formula for 
the character values of simple $\B'_k$-modules 
, in a form adapted to the simple modules in the ${\Bbb Z}_2$-graded sense. 
\proclaim{Theorem 4.5}{\rm\!\!(cf\. [10, Lem\. 7.5])} 
We have 
$$
\align
&2^{\frac{3(l(\kappa)+l(\nu))}{2}}
p_{\kappa}(x)p_{\nu}(y)\\
&\quad\quad=
\sum_{(\lam,\mu)\in{(DP^2)}_k}
\Ch[V_{\lam,\mu}](w^{\kappa,\nu})
2^{\frac{-l(\lam)-l(\mu)-\vep(\lam,\mu)}{2}}
Q_{\lam}(x,y)Q_{\mu}(x,-y)
\endalign
$$
for all $(\kappa,\nu)\in{(OP^2)}_k$, where 
$Q_{\lam}(x,y)=Q_{\lam}(x_1,x_2,\dots,y_1,y_2,\dots)$ and 
$Q_{\mu}(x,-y)$
$=Q_{\mu}(x_1,x_2,\dots,-y_1,-y_2,\dots)$. 
\endproclaim
The formula \thetag{4.15} is different from Stembridge's formula.
Let us mention a relationship between the two formulas. 
\comment
Define an automorphism $f(x,y)\mapsto f^{*}(x,y)$ of 
$\Omega_x\otimes \Omega_y$ by 
$p^{*}_r(x)=p_r(x)$ and $p^{*}_r(y)=-p_r(y)$ for all odd $r$ 
(in the notation of [6], we have $f^{*}(y)=\omega f(-y)$). 
We have $f^{*}(x,y)=f(x,-y)$ for $f(x,y)\in\Omega_x\otimes\Omega_y$. 
Moreover, 
\endcomment
Define an algebra endomorphism $\iota$ of 
$\Omega_x\otimes\Omega_y$ by 
$\iota(f\otimes 1)=f(x,y)=f(x_1,x_2,\dots,y_1,y_2,\dots)$ and 
$\iota(1\otimes g)=g(x,-y)=g(x_1,x_2,\dots,-y_1,-y_2,\dots)$. 
\comment
 (since the 
$y$-part belongs to $\Omega_y$, this ``na\"\i ve notation'' 
actually coincides with $g(x-y)$ in the $\Lam$-ring notation). 
\endcomment
Note that 
$\{Q_{\lam}(x,y)Q_{\mu}(x,-y)\;|\;(\lam,\mu)\in{(DP^2)}\}$ is 
a basis of $\Omega_x\otimes\Omega_y$ (cf\. [10, Th\. 7.1, Lem\. 7.5]). 
It follows that $\iota$ is an automorphism, 
since $\iota(Q_{\lam}(x)Q_{\mu}(y))=Q_{\lam}(x,y)Q_{\mu}(x,-y)$. 
Moreover, since $\iota(p_r(x,y))=2p_r(x)$ and $\iota(p_r(x,-y))=2p_r(y)$ 
for any odd $r$, it follows that 
the image of \thetag{4.15} under $\iota$ coincides with 
Stembridge's formula. 

\head
Appendix
\endhead
\subhead
A. Sergeev's duality
\endsubhead 
We review Sergeev's duality relation between $\B'_k$ and $\U_n$ 
using DSCT. Define a map $d\:DP_k\to{\Bbb Z}_2$ by 
$d(\lam)=0$ (resp\. $d(\lam)=1$) 
if $l(\lam)$ is even (resp\. $l(\lam)$ is odd). 

\proclaim{Theorem A}{\rm\!\![8]} 
\t{\rm(1)} The two superalgebras 
$\Psi(\B_k)$ and $\U_n$ act on $W$ as mutual supercentralizers 
of each other{\rm:} 
$$
\E^{\cdot}_{\Theta(\U_n)}(W)=\Psi(\B_k),
\quad \E^{\cdot}_{\Psi(\B_k)}(W)=\Theta(\U_n).\tag A.1
$$

\t{\rm(2)} The simple $\B_k$-module $W_{\lam}$ 
$(\lam\in DP_k)$ occurs in $W$ if and only if $l(\lam))\leq n$.
 Then we have 
$$
W\cong_{\B_k\lot \U_n}
\bigoplus_{
\smallmatrix
\lam \in DP_k\\
l(\lam)\leq n
\endsmallmatrix}
W_{\lam}\hotimes U_{\lam}\tag A.2
$$
where $U_{\lam}$ denotes 
a simple $\U_n$-module corresponding to $W_{\lam}$ in $W$
 in the sense of DSCT. 

\t{\rm(3)} The character values of $\Ch[U_{\lam}]$ are given as follows{\rm:} 
$$
\Ch[U_{\lam}](x_1,x_2,\dots,x_n)
=(\sqrt{2})^{d(\lam)-l(\lam)}Q_{\lam}(x_1,x_2,\dots,x_n).\tag A.3
$$
\endproclaim 
\subhead
B. A duality of $\A_k$ and $\fq(n)$.
\endsubhead
We established a duality relation between 
$\A_k$ and $\U_n$ on the same space $W^{\vep}$ as in Theorem 4.2. 
\proclaim{Theorem B}{\rm\!\! [11, Th\. 4.1]} 
The submodule $W^{\vep}$ is decomposed 
as a multiplicity-free sum of simple 
$\A_k\lot\U_n$-modules as follows{\rm:} 
$$
W^{\vep}\cong_{\A_k\lot\U_n}
\bigoplus_{\lam\in DP_k}
V_{\lam}\hotimes U_{\lam}.
\tag B.1
$$

\t{\rm(1)} Assume that $k$ is even. Then the simple $\A_k\lot\U_n$-modules 
$V_{\lam}\hotimes U_{\lam}$ in $W^{\vep}$ are of type $M$. 
Furthermore we have 
$$
\E^{\cdot}_{\Theta(\U_n)}(W^{\vep})=\Psi(\A_k),\quad\;
\E^{\cdot}_{\Psi(\A_k)}(W^{\vep})=\Theta(\U_n).\tag B.2
$$

\t{\rm(2)} Assume that $k$ is odd. Then the simple $\A_k\lot\U_n$-modules 
$V_{\lam}\hotimes U_{\lam}$ in $W^{\vep}$ are of type $Q$. 
Furthermore we have 
$$
\E^{\cdot}_{\Theta(U_n)}(W^{\vep})\cong\C_1\otimes\Psi(\A_k),\quad
\E^{\cdot}_{\Psi(\A_k)}(W^{\vep})\cong\C_1\otimes\Theta(\U_n).
\tag B.3
$$
\endproclaim

\enddocument